\documentclass[11pt]{article}
\usepackage[a4paper,margin=28mm]{geometry}
\usepackage{iftex}

\ifPDFTeX
  \usepackage[T1]{fontenc}
  \usepackage{libertinus}
  \usepackage{amsmath,amssymb,amsthm,mathtools}
  \usepackage{mathrsfs}
\else
  \usepackage{libertinus-otf}
  \usepackage{amsthm,mathtools}
\fi

\usepackage{microtype}
\usepackage{enumitem}
\usepackage{booktabs}
\usepackage{hyperref}
\usepackage{titlesec}
\usepackage{indentfirst}
\allowdisplaybreaks

\hypersetup{
  hidelinks,
  pdftitle={The Capillary Orlicz-Christoffel-Minkowski Problem},
  pdfsubject={Research article},
  pdfauthor={Hua Zhu},
  pdfkeywords={capillary hypersurfaces; Orlicz-Christoffel-Minkowski problem; k-th area measure; fully nonlinear elliptic equations; continuity method; constant-rank theorem}
}

\titleformat{\section}{\large\bfseries}{\thesection}{0.75em}{}
\titleformat{\subsection}{\normalsize\bfseries}{\thesubsection}{0.7em}{}
\titleformat{\subsubsection}{\normalsize\itshape}{\thesubsubsection}{0.7em}{}
\titlespacing*{\section}{0pt}{2.7ex plus .8ex minus .2ex}{1.1ex}
\titlespacing*{\subsection}{0pt}{2.1ex plus .6ex minus .2ex}{.8ex}
\titlespacing*{\subsubsection}{0pt}{1.7ex plus .5ex minus .2ex}{.6ex}

\newtheoremstyle{submissionplain}
  {6pt}{6pt}
  {\itshape}{}
  {\bfseries}{.}
  {.5em}{}
\newtheoremstyle{submissiondef}
  {6pt}{6pt}
  {}{}
  {\bfseries}{.}
  {.5em}{}
\newtheoremstyle{submissionremark}
  {5pt}{5pt}
  {}{}
  {\itshape}{.}
  {.5em}{}

\theoremstyle{submissionplain}
\newtheorem{theorem}{Theorem}[section]
\newtheorem{proposition}[theorem]{Proposition}
\newtheorem{lemma}[theorem]{Lemma}
\newtheorem{corollary}[theorem]{Corollary}
\theoremstyle{submissiondef}

\theoremstyle{submissionremark}
\newtheorem{remark}[theorem]{Remark}

\newcommand{\Ct}{\mathcal C_{\theta}}
\newcommand{\sig}{\sigma}
\newcommand{\W}{W}
\newcommand{\ellc}{\ell}
\newcommand{\bk}{\bar\sigma_k}

\newcommand{\R}{\mathbb R}
\newcommand{\dd}{\,d}

\newcommand{\rank}{\operatorname{rank}}
\newcommand{\diag}{\operatorname{diag}}
\newcommand{\eps}{\varepsilon}

\begin{document}

\begin{center}
{\LARGE\bfseries The Capillary Orlicz--Christoffel--Minkowski Problem}\par\vspace{5mm}
{\large Hua Zhu\textsuperscript{1,2}}\par\vspace{3mm}
{\small
\textsuperscript{1}Applied Mathematics Center of China (Mianyang) Science and Technology City, Mianyang, Sichuan, \\
621010, People's Republic of China.\\
\par\vspace{1.5mm}
\textsuperscript{2}School of Mathematics and Physics, Southwest University of Science and Technology, Mianyang, Sichuan,\\
621010, People's Republic of China.\\
Email: \emph{zhuhmaths@mail.ustc.edu.cn}}\par
\end{center}
\vspace{5mm}

\begin{abstract}
We study the Capillary Orlicz--Christoffel--Minkowski Problem on an acute spherical cap. For $1\le k\le n-1$, the equation is
\[
h\,\varphi(\ell/h)\,\bar\sigma_k(\nabla^2h+h\sigma)=f,
\qquad \nabla_\mu h=\cot\theta\,h,
\]
where $\ell$ is the support function of the unit capillary spherical cap. For $1\le k\le n-1$, this includes the capillary Christoffel--Minkowski and capillary $L_p$ Christoffel--Minkowski problems; the same density formally recovers the capillary Orlicz--Minkowski endpoint when $k=n$. We organize the theory in five parts. After recalling the elementary symmetric functions and the capillary support-function formalism, we review the capillary Orlicz--Brunn--Minkowski framework and introduce the intermediate-order Orlicz $k$-mixed functional and its Euler--Lagrange equation. In the supercritical regime
\[
a_\varphi(t):=\frac{t\varphi'(t)}{\varphi(t)}\ge k+1+\delta,
\]
we derive uniform $C^0$, $C^1$ and $C^2$ estimates for positive $k$-admissible solutions. The same supercritical index produces a strictly coercive second variation and a spectral gap for the linearized Neumann operator. These facts lead to existence and uniqueness by a continuity method, quantitative stability of the solution map in both the datum and the Orlicz weight, and an open perturbative class of strictly convex non-power solutions. Under an additional capillary Orlicz-convexity condition on the datum, a full-rank argument upgrades the global admissible solution to a strictly convex capillary support function.
\end{abstract}


\section{Introduction}

Prescribed curvature and area-measure problems form a central part of convex geometry and fully nonlinear geometric analysis. In the closed setting, the Minkowski and Christoffel--Minkowski problems, together with their $L_p$ and Orlicz extensions, can be expressed through nonlinear equations for the support function. The $L_p$ Brunn--Minkowski theory originates in the Firey-type framework developed by Lutwak \cite{Lutwak93}, while the Orlicz Minkowski problem and the broader Orlicz--Brunn--Minkowski theory were developed in \cite{HLYZ,GHW,XiJinLeng}. For the Christoffel--Minkowski problem, the Hessian-equation approach and the convexity and existence theory were developed in particular by Guan--Ma and Guan--Ma--Zhou \cite{GuanMa,GuanMaZhou}; the general $k$-Hessian elliptic framework goes back to Caffarelli--Nirenberg--Spruck \cite{CNS}. For hypersurfaces with capillary boundary, the same geometric quantities are coupled to the contact-angle condition, and the corresponding support-function equations acquire a nontrivial boundary condition. This boundary interaction is the main analytic feature distinguishing the capillary problem from its closed counterpart. Standard regularity theory for nonlinear oblique boundary problems provides the analytic background for such boundary conditions \cite{LiebermanTrudinger}.

Let
\[
\R^{n+1}_+=\{x\in\R^{n+1}:x_{n+1}>0\},
\qquad e=-E_{n+1},
\]
and let $\Sigma\subset\overline{\R^{n+1}_+}$ be a smooth strictly convex capillary hypersurface meeting $\partial\R^{n+1}_+$ at a constant angle $\theta\in(0,\pi)$. The basic model is the unit capillary spherical cap $\Ct$. Through the capillary Gauss map, a strictly convex capillary hypersurface can be parametrized over $\Ct$. We denote its capillary support function by $h$ and set
\[
\ellc(\xi)=\sin^2\theta+\cos\theta\langle\xi,e\rangle,
\qquad
\W[h]=\nabla^2h+h\sig.
\]
The eigenvalues of $\W[h]$ are the principal radii of $\Sigma$, while the capillary condition becomes the Robin relation
\[
\nabla_\mu h=\cot\theta\,h
\qquad\text{on }\partial\Ct.
\]
The capillary $k$-th area measure has density proportional to
\[
\ellc\,\bk(\W[h])\dd\sig,
\qquad
\bk(A):=\binom nk^{-1}\sigma_k(A).
\]

The recent capillary Brunn--Minkowski theory also includes Minkowski-type inequalities for capillary hypersurfaces \cite{WangWengXia} and, in a complementary subcritical range, the even capillary $L_p$-curvature problem \cite{HuIvaki}. Motivated by the capillary Orlicz--Minkowski theory at the Monge--Amp\`ere endpoint $k=n$ \cite{WangZhu} and by the capillary $L_p$ Christoffel--Minkowski theory \cite{MWWLp,CMCap}, we introduce an Orlicz modification of this intermediate-order area measure.

\medskip
\noindent\textbf{Capillary Orlicz $k$-th area measure.}
For $1\le k\le n-1$ and a positive Orlicz weight $\varphi$, define
\begin{equation}\label{intro-measure}
\dd S^c_{\varphi,k}[h]
:=h\,\varphi\!\left(\frac{\ellc}{h}\right)
\bk(\W[h])\dd\sig.
\end{equation}
This definition contains the familiar cases already present in the capillary theory. Indeed,
\begin{align*}
\varphi(t)=t&\quad\Longrightarrow\quad
\dd S^c_{\varphi,k}=\ellc\,\bk(\W[h])\dd\sig,\\
\varphi(t)=t^p&\quad\Longrightarrow\quad
\dd S^c_{\varphi,k}=\ellc^p h^{1-p}\bk(\W[h])\dd\sig,\\
k=n&\quad\Longrightarrow\quad
\text{the capillary Orlicz--Minkowski density of \cite{WangZhu}.}
\end{align*}
In particular, when $\varphi(t)=t^p$ and $k=n$, since $\bk(W)=\det W$, the measure above becomes
\[
\dd S^c_{\varphi,n}
=\ellc^p h^{1-p}\det(\nabla^2h+h\sig)\dd\sig,
\]
which is exactly the capillary $L_p$-surface area measure introduced by Mei--Wang--Weng \cite{MWWLpMinkowski}. Their prescribed problem is written as $\dd S_p^c=\widetilde f\,\ellc^p\dd\sig$. Thus, after the identification $f=\ellc^p\widetilde f$, the prescription \eqref{intro-prescription} reduces precisely to
\[
\det(\nabla^2h+h\sig)=\widetilde f\,h^{p-1},
\qquad
\nabla_\mu h=\cot\theta\,h,
\]
the capillary $L_p$-Minkowski equation of \cite{MWWLpMinkowski}. For $1\le k\le n-1$, the same power choice $\varphi(t)=t^p$ gives the capillary $L_p$ Christoffel--Minkowski measure. For $k=n$, Wang--Zhu developed the capillary Orlicz--Minkowski and Orlicz--Brunn--Minkowski theory \cite{WangZhu}. For intermediate orders, Mei--Wang--Weng developed the corresponding supercritical capillary $L_p$ curvature theory and the estimates used for the $L_p$ Christoffel--Minkowski equation \cite{MWWLp}; the capillary Christoffel--Minkowski problem was studied in \cite{CMCap}. Hu--Ivaki--Scheuer obtained a complementary capillary Christoffel--Minkowski existence result under a spherical-convexity condition \cite{HuIvakiScheuer}, while Hu--Ivaki treated the even capillary $L_p$ Christoffel--Minkowski problem in the subcritical range $1<p<k+1$ \cite{HuIvakiLpCM}. These results complement the supercritical regime considered below.

This leads naturally to the Capillary Orlicz--Christoffel--Minkowski Problem considered here.

\medskip
\noindent\textbf{Capillary Orlicz--Christoffel--Minkowski Problem.}
Given a positive smooth function $f$ on $\overline\Ct$, can one find a smooth strictly convex capillary hypersurface whose capillary Orlicz $k$-th area measure satisfies
\begin{equation}\label{intro-prescription}
\dd S^c_{\varphi,k}[h]=f\dd\sig?
\end{equation}
Equivalently, the problem is to determine a capillary support function $h$ realizing the prescribed density $f$.

In view of the support-function representation recalled above, the geometric prescription has a direct analytic formulation. Combining \eqref{intro-measure} with \eqref{intro-prescription} gives
\[
h\,\varphi(\ellc/h)\bk(\nabla^2h+h\sig)=f
\qquad\text{in }\Ct,
\]
while capillarity imposes the Robin boundary condition. Consequently the Capillary Orlicz--Christoffel--Minkowski Problem is equivalent to the following $k$-Hessian-type elliptic boundary value problem:
\begin{equation}\label{main-eq}
\begin{cases}
\displaystyle
\bk(\nabla^2h+h\sig)
=\dfrac{f}{h\,\varphi(\ellc/h)},&\text{in }\Ct,\\[2mm]
\nabla_\mu h=\cot\theta\,h,&\text{on }\partial\Ct.
\end{cases}
\end{equation}
We call a solution $k$-admissible if
\[
\lambda(\nabla^2h+h\sig)\in\Gamma_k.
\]
Convexity and strict convexity correspond respectively to $\nabla^2h+h\sig\ge0$ and $\nabla^2h+h\sig>0$. Thus the prescribed geometric measure problem is reduced to a fully nonlinear $k$-Hessian-type equation on the spherical cap with a Robin boundary condition.

We now describe the parameter range considered in this paper. Define the logarithmic Orlicz slope by
\begin{equation*}
a_\varphi(t):=\frac{t\varphi'(t)}{\varphi(t)}.
\end{equation*}
We work in the supercritical regime
\begin{equation}\label{SC}
a_\varphi(t)\ge k+1+\delta\qquad\forall t>0
\end{equation}
for some $\delta>0$. When $\varphi(t)=t^p$, this reduces exactly to $p>k+1$, the supercritical range of the capillary $L_p$ Christoffel--Minkowski equation. The same quantity also appears in the basic dilation law. If $h_c=ch$ and $r=\ellc/h$, then
\begin{equation*}
\frac{h_c\varphi(\ellc/h_c)\bk(\W[h_c])}
{h\varphi(\ellc/h)\bk(\W[h])}
=c^{k+1}\frac{\varphi(r/c)}{\varphi(r)}.
\end{equation*}
Accordingly, it is convenient to set
\begin{equation*}
\Lambda_\varphi(t):=a_\varphi(t)-(k+1).
\end{equation*}
Then
\[
\left.\frac{d}{d\log c}\right|_{c=1}
\log\!\left[c^{k+1}\varphi(r/c)\right]
=-\Lambda_\varphi(r).
\]
Hence \eqref{SC} expresses strict monotonicity of the prescribed density under dilation. The identity $\Lambda_\varphi\equiv0$ gives $\varphi(t)=Ct^{k+1}$, while logarithmic interpolation preserves a fixed supercritical margin: if $\varphi_s=\varphi_0^{1-s}\varphi_1^s$, then
\[
a_{\varphi_s}=(1-s)a_{\varphi_0}+sa_{\varphi_1}.
\]

We now state the main result. For the strict-convexity conclusion, define
\[
H:=f\ellc^{k-1}.
\]
We shall use the structural conditions
\begin{equation}\label{OC-intro}
\nabla^2(\ellc H^{-\alpha})+\ellc H^{-\alpha}\sig\ge0
\qquad\text{for every }0<\alpha\le\frac1k,
\end{equation}
and
\begin{equation}\label{OC-boundary-intro}
H_\mu\ge0\qquad\text{on }\partial\Ct.
\end{equation}

The main theorem of this paper is as follows.

\begin{theorem}\label{thm:mainA}
Let $0<\theta<\pi/2$ and $1\le k\le n-1$. Suppose that $f\in C^\infty(\overline\Ct)$ is positive and $\varphi\in C^\infty((0,\infty))$ is positive. Assume that \eqref{SC} holds,
\[
t a_\varphi'(t)\ge0\qquad(t>0),
\]
and that $H=f\ellc^{k-1}$ satisfies \eqref{OC-intro}--\eqref{OC-boundary-intro}. Then there exists a smooth strictly convex capillary hypersurface $\Sigma\subset\overline{\R^{n+1}_+}$ with contact angle $\theta$ such that its capillary support function $h$ solves the Capillary Orlicz--Christoffel--Minkowski Problem, namely
\[
\dd S^c_{\varphi,k}[h]=f\dd\sig
\qquad\text{on }\Ct.
\]
Moreover, $\Sigma$ is unique.
\end{theorem}

\begin{theorem}[Perturbative stability of strict convexity]\label{thm:mainC}
Let $(f_0,\varphi_0)$ be positive smooth data satisfying \eqref{SC} with a fixed supercritical margin $\delta>0$, and suppose the corresponding admissible solution $h_0$ is strictly convex. Fix $0<\alpha<1$ and a compact interval $I\Subset(0,\infty)$ whose interior contains the range of $\ellc/h_0$. Then there exists $\varepsilon>0$ such that every positive smooth pair $(f,\varphi)$ satisfying the same supercritical margin and 
\[
\|f-f_0\|_{C^\alpha(\overline\Ct)}
+
\|\varphi-\varphi_0\|_{C^2(I)}<\varepsilon
\]
has a unique smooth strictly convex solution $h$ with $\ellc/h\in I$. In particular, every strictly convex supercritical power-law solution has a neighborhood of genuinely non-power Orlicz solutions.
\end{theorem}

The proof follows an elliptic continuity scheme. Section 2 collects the elementary symmetric-function identities and the capillary support-function geometry. Section 3 recalls the capillary Orlicz--Brunn--Minkowski framework and the capillary mixed-volume formalism, and introduces the intermediate-order Orlicz mixed functional and its variational structure. Section 4 establishes the uniform $C^0$, $C^1$, and $C^2$ estimates. Section 5 proves the main results by the continuity method, establishes the coercivity and stability properties of the linearized problem, and then uses the full-rank continuation to obtain global strict convexity.

\begin{remark}[Intermediate-order Orlicz geometry]
The $k=n$ Orlicz combination and Orlicz--Brunn--Minkowski inequalities quoted in Section 3 are established results of Wang--Zhu \cite{WangZhu}. For $1\le k\le n-1$, this manuscript uses the same first-variation pattern to define the Orlicz $k$-mixed functional and the measure \eqref{intro-measure}. A complete intermediate-order capillary Orlicz--Brunn--Minkowski inequality is not needed for Theorem \ref{thm:mainA} and is not claimed here.
\end{remark}

\section{Preliminaries}

\subsection{Elementary symmetric functions and the G\aa rding cones}

For $\lambda=(\lambda_1,\dots,\lambda_n)\in\R^n$, let
\[
\sigma_m(\lambda)
=\sum_{1\le i_1<\cdots<i_m\le n}
\lambda_{i_1}\cdots\lambda_{i_m},
\qquad 0\le m\le n,
\]
with $\sigma_0=1$. If $A$ is a symmetric endomorphism, $\sigma_m(A)$ denotes the same polynomial of its eigenvalues. We also use
\[
\bar\sigma_m(A):=\binom nm^{-1}\sigma_m(A).
\]
The $k$-th G\aa rding cone is
\[
\Gamma_k:=\{\lambda\in\R^n:\sigma_j(\lambda)>0,\ 1\le j\le k\},
\]
where the connected component containing the positive cone is understood. A symmetric tensor $A$ is called $k$-admissible if $\lambda(A)\in\Gamma_k$.

At a diagonal matrix $A=\diag(\lambda_1,\dots,\lambda_n)$,
\begin{equation*}
\frac{\partial\sigma_k}{\partial A_{ii}}=\sigma_{k-1}(\lambda|i),
\qquad
\frac{\partial\sigma_k}{\partial A_{ij}}=0\quad(i\ne j),
\end{equation*}
where $\sigma_{k-1}(\lambda|i)$ omits $\lambda_i$. Euler homogeneity gives
\begin{equation*}
\sigma_k^{ij}(A)A_{ij}=k\sigma_k(A).
\end{equation*}
In particular, if
\[
F(A)=\bar\sigma_k(A)^{1/k},
\qquad
\mathcal F(A)=\log\bar\sigma_k(A),
\]
then
\begin{equation*}
F^{ij}A_{ij}=F,
\qquad
\mathcal F^{ij}A_{ij}=k.
\end{equation*}

We shall repeatedly use the following standard properties.

\begin{lemma}[Newton--Maclaurin and concavity]\label{lem:symmetric}
Let $\lambda\in\Gamma_k$.
\begin{enumerate}[label=(\roman*),leftmargin=2.0em]
\item If $1\le j<k$, then
\[
\bar\sigma_j(\lambda)^{1/j}\ge \bar\sigma_k(\lambda)^{1/k}.
\]
\item $\sigma_k^{ij}$ is positive definite as a linearized coefficient matrix on $\Gamma_k$.
\item $F=\bar\sigma_k^{1/k}$ is concave and one-homogeneous on $\Gamma_k$.
\item $F$ is inverse-concave in the sense used in the constant-rank theory of \cite{BISCRT}.
\end{enumerate}
\end{lemma}

\begin{lemma}[Compactness inside the G\aa rding cone]\label{lem:Gamma-compact}
Fix $c_*>0$ and $C_*<\infty$. The set
\[
\{\lambda\in\Gamma_k: |\lambda|\le C_*,\ \bar\sigma_k(\lambda)\ge c_*\}
\]
is contained in a compact subset of $\Gamma_k$. Consequently $\bar\sigma_k^{1/k}$ is uniformly elliptic on this set.
\end{lemma}

\begin{proof}
The closure of the displayed set inside the Euclidean ball $\{|\lambda|\le C_*\}$ is compact. Since $\Gamma_k$ is the component of $\{\sigma_k>0\}$ containing the positive cone, its finite boundary is contained in $\{\sigma_k=0\}$. The lower bound on $\bar\sigma_k$ therefore prevents the closure from meeting $\partial\Gamma_k$. Uniform ellipticity follows by continuity of the positive matrix $F^{ij}$ on the resulting compact subset.
\end{proof}

For later use, denote the Newton tensor
\[
S_k^{ij}(A):=\frac{\partial\bar\sigma_k}{\partial A_{ij}}.
\]
If $A=\nabla^2h+h\sig$ on the unit sphere, then the Codazzi identities imply
\begin{equation}\label{Newton-divfree}
\nabla_iS_k^{ij}(A)=0.
\end{equation}
This divergence-free property is the basic reason that mixed-volume expressions are symmetric after integration by parts.

\subsection{Capillary hypersurfaces and the spherical cap}

Let
\[
\Ct:=\{\xi\in\overline{\R^{n+1}_+}:|\xi-\cos\theta\,e|=1\},
\qquad e=-E_{n+1}.
\]
Its round metric is denoted by $\sig$, and $\mu$ is the outward unit conormal to $\partial\Ct\subset\Ct$. The support function of the unit capillary spherical cap is
\begin{equation*}
\ellc(\xi)=\sin^2\theta+\cos\theta\langle\xi,e\rangle.
\end{equation*}
It satisfies
\begin{equation}\label{ell-identities}
\nabla^2\ellc+\ellc\sig=\sig,
\qquad
\nabla_\mu\ellc=\cot\theta\,\ellc.
\end{equation}

Let $\widehat\Sigma$ be a smooth strictly convex capillary body, and let $\nu$ be the outer unit normal on the curved part $\Sigma$. The capillary Gauss map is
\[
\widetilde\nu=\nu+\cos\theta\,e:\Sigma\longrightarrow\Ct.
\]
For a strictly convex capillary hypersurface it is a diffeomorphism. Writing $X=\widetilde\nu^{-1}(\xi)$, define
\begin{equation*}
h(\xi)=\langle X(\xi),\xi-\cos\theta e\rangle.
\end{equation*}
Then \cite{MinkowskiCap,AFII}
\begin{equation*}
X(\xi)=\nabla h+h(\xi-\cos\theta e),
\end{equation*}
\begin{equation}\label{Robin}
\nabla_\mu h=\cot\theta\,h\quad\text{on }\partial\Ct,
\end{equation}
and the eigenvalues of
\begin{equation*}
\W[h]:=\nabla^2h+h\sig
\end{equation*}
are exactly the principal radii of $\Sigma$. Hence strict convexity is equivalent to $\W[h]>0$.

It is convenient to divide by the spherical-cap support function:
\begin{equation*}
u:=\frac{h}{\ellc}.
\end{equation*}
By \eqref{ell-identities} and \eqref{Robin},
\begin{equation}\label{u-Neumann}
\nabla_\mu u=0\quad\text{on }\partial\Ct.
\end{equation}
Moreover,
\begin{equation}\label{Btheta}
\W[h]=\mathcal B_\theta[u]
:=\ellc\nabla^2u+\nabla\ellc\otimes\nabla u+\nabla u\otimes\nabla\ellc+u\sig.
\end{equation}
This is the Neumann formulation used in the $C^0$ estimate.

\subsection{Boundary identities}

Let $\{e_\alpha\}_{\alpha=1}^{n-1}$ be an orthonormal frame tangent to $\partial\Ct$ and set $e_n=\mu$. Since $\partial\Ct$ is totally umbilical in $\Ct$ with principal curvature $\cot\theta$, any capillary function satisfying \eqref{Robin} obeys
\begin{equation*}
\nabla^2h(e_\alpha,\mu)=0.
\end{equation*}
Thus
\begin{equation}\label{W-mixed-zero}
\W_{\alpha n}=0\qquad\text{on }\partial\Ct.
\end{equation}
Differentiating once more yields the standard relation
\begin{equation}\label{W-normal-derivative}
\W_{\alpha\alpha n}
=\cot\theta\big(\W_{nn}-\W_{\alpha\alpha}\big)
\qquad\text{on }\partial\Ct,
\end{equation}
which is fundamental both in the boundary $C^2$ estimate and in the full-rank argument.

\subsection{Mixed discriminants and capillary mixed volumes}

For symmetric $n\times n$ matrices $A_1,\dots,A_m$, their mixed discriminant $Q$ is characterized by
\[
\det(\lambda_1A_1+\cdots+\lambda_mA_m)
=\sum_{i_1,\dots,i_n=1}^{m}\lambda_{i_1}\cdots\lambda_{i_n}
Q(A_{i_1},\dots,A_{i_n}).
\]
For capillary functions $f_1,\dots,f_{n+1}$, define
\begin{equation*}
V(f_1,\dots,f_{n+1})
:=\frac1{n+1}\int_{\Ct}
f_1Q(\W[f_2],\dots,\W[f_{n+1}])\dd\sig.
\end{equation*}
The divergence-free identity for the mixed Newton tensor together with the common Robin boundary condition makes $V$ symmetric and multilinear; see \cite{AFII}. For support functions it agrees with the mixed volume of capillary convex bodies. In particular, the classical capillary Alexandrov--Fenchel inequality holds for all contact angles $\theta\in(0,\pi)$ \cite{AFII}.

We shall also use the capillary Minkowski identity: if $h$ is a capillary support function, then for $0\le r\le n-1$,
\begin{equation}\label{capillary-Minkowski-identity}
(n-r)\int_{\Ct}h\,\sigma_r(\W[h])\dd\sig
=(r+1)\int_{\Ct}\ellc\,\sigma_{r+1}(\W[h])\dd\sig.
\end{equation}
It will exclude a non-full constant rank in Section 5.

\section{Capillary Orlicz--Brunn--Minkowski Theory}

\subsection{Capillary Orlicz combination}

We first recall the geometric Orlicz addition used by Wang--Zhu \cite{WangZhu}, itself modeled on the Orlicz--Brunn--Minkowski theory of Gardner--Hug--Weil \cite{GHW}. In this subsection let $\Phi:[0,\infty)\to[0,\infty)$ be convex, strictly increasing, $C^1$ on $(0,\infty)$, and normalized by $\Phi(0)=0$ and $\Phi(1)=1$. The normalization is only for the geometric addition; the PDE weight $\varphi$ used in \eqref{main-eq} need not be normalized.

Let $K,L$ be capillary convex bodies whose support functions $h_K,h_L$ are positive on $\Ct$, and let $\alpha,\beta\ge0$ with $\alpha+\beta>0$. Define $M_\Phi(\alpha,\beta;K,L)$ by its support function
\begin{equation}\label{Orlicz-combination}
h_{M_\Phi}(\xi)
:=\inf\left\{s>0:
\alpha\Phi\!\left(\frac{h_K(\xi)}s\right)
+\beta\Phi\!\left(\frac{h_L(\xi)}s\right)\le1\right\}.
\end{equation}
Equivalently, wherever $h=h_{M_\Phi}>0$,
\begin{equation}\label{Orlicz-implicit}
\alpha\Phi(h_K/h)+\beta\Phi(h_L/h)=1.
\end{equation}

\begin{proposition}[Robin closure of the Orlicz combination]\label{prop:Orlicz-Robin}
If $h_K$ and $h_L$ satisfy the capillary Robin condition, then the support function defined by \eqref{Orlicz-combination} also satisfies
\[
\nabla_\mu h=\cot\theta\,h
\qquad\text{on }\partial\Ct.
\]
Consequently the capillary class is closed under the Orlicz combination whenever the support function in \eqref{Orlicz-combination} is spherically convex.
\end{proposition}

\begin{proof}
Differentiate \eqref{Orlicz-implicit} in the $\mu$ direction. Since
\[
\left(\frac{h_K}{h}\right)_\mu
=\frac{h_K}{h^2}(\cot\theta h-h_\mu),
\qquad
\left(\frac{h_L}{h}\right)_\mu
=\frac{h_L}{h^2}(\cot\theta h-h_\mu),
\]
we obtain
\[
\left[
\alpha h_K\Phi'(h_K/h)+\beta h_L\Phi'(h_L/h)
\right](\cot\theta h-h_\mu)=0.
\]
The bracket is positive because $\Phi$ is strictly increasing. Hence $h_\mu=\cot\theta h$. The spherical convexity and reconstruction of a capillary body then follow from the support-function correspondence in \cite{AFII,WangZhu}.
\end{proof}

When $\Phi(t)=t^p$, \eqref{Orlicz-combination} reduces to the capillary $L_p$ combination. Thus the Robin condition is compatible not only with Minkowski addition but with the nonlinear Orlicz addition as well. This is the geometric reason that a capillary Orlicz theory is natural.

\subsection{The known capillary Orlicz endpoint}

For $k=n$, Wang--Zhu define a capillary Orlicz mixed volume through the first variation of volume under the Orlicz perturbation. With a harmless choice of normalization, the smooth mixed functional can be written as
\begin{equation}\label{WZ-Orlicz-volume}
V_\Phi^c(K,L)
:=\frac1{n+1}\int_{\Ct}
\Phi\!\left(\frac{h_L}{h_K}\right)
h_K\det\W[h_K]\dd\sig.
\end{equation}
Their first-variation formula contains the conventional factor $1/\Phi'(1)$, which we suppress in \eqref{WZ-Orlicz-volume}. Combining the capillary Alexandrov--Fenchel inequality with Jensen's inequality, they prove capillary Orlicz--Minkowski and Orlicz--Brunn--Minkowski inequalities and characterize the equality case under the corresponding strict-convexity hypotheses \cite{WangZhu}. These results form the $k=n$ geometric endpoint of the present program.

The full intermediate-order analogue of those inequalities is a separate geometric problem. For the elliptic existence theorem we only need the first-variation pattern suggested by \eqref{WZ-Orlicz-volume}, together with the capillary $k$-th area measure of \cite{CMCap,MWWLp}.

\subsection{An intermediate-order Orlicz mixed functional}

For a positive capillary support function $h_K$ and another positive support function $h_L$, define
\begin{equation*}
V^c_{\Phi,k}(K,L)
:=\frac1{n+1}\int_{\Ct}
\Phi\!\left(\frac{h_L}{h_K}\right)
h_K\bk(\W[h_K])\dd\sig.
\end{equation*}
This normalization is chosen so that:
\begin{itemize}[leftmargin=1.8em]
\item for $k=n$, it coincides with \eqref{WZ-Orlicz-volume} up to the normalization of $\sigma_n$;
\item for $\Phi(t)=t$, the integrand becomes $h_L\bk(\W[h_K])$;
\item taking $L$ to be the unit spherical cap, $h_L=\ellc$, gives precisely the Orlicz factor appearing in the proposed measure \eqref{intro-measure}.
\end{itemize}

Accordingly, for the PDE weight $\varphi$ we define
\begin{equation}\label{Orlicz-k-measure}
\dd S^c_{\varphi,k}[h]
=h\,\varphi(\ellc/h)\bk(\W[h])\dd\sig.
\end{equation}
The equation \eqref{main-eq} is therefore the prescribed density problem for \eqref{Orlicz-k-measure}.

\begin{remark}
Equation \eqref{Orlicz-k-measure} is the natural common extension of two established measures: the capillary $k$-th area measure $\ellc\,\sigma_k(\W[h])\dd\sig$ and the capillary Orlicz surface-area measure $h\varphi(\ellc/h)\det\W[h]\dd\sig$. The normalization by $\binom nk$ changes no solvability issue.
\end{remark}

\subsection{Variational structure}

The mixed-volume symmetry from Section 2 produces a useful functional. Define
\begin{equation}\label{Mk-functional}
\mathcal M_k[h]
:=\frac1{k+1}\int_{\Ct}h\bk(\W[h])\dd\sig.
\end{equation}

\begin{lemma}[First variation of $\mathcal M_k$]\label{lem:first-var-M}
If $\eta$ satisfies the same Robin boundary condition as $h$, then
\begin{equation}\label{first-var-M}
D\mathcal M_k[h](\eta)
=\int_{\Ct}\eta\bk(\W[h])\dd\sig.
\end{equation}
\end{lemma}

\begin{proof}
Differentiating \eqref{Mk-functional} gives
\[
\frac1{k+1}\int
\left[\eta\bk(\W)+hS_k^{ij}(\eta_{ij}+\eta\sig_{ij})\right]\dd\sig.
\]
By the divergence-free identity \eqref{Newton-divfree}, two integrations by parts, and the common Robin condition, the second term equals
\[
\frac1{k+1}\int \eta S_k^{ij}\W_{ij}\dd\sig
=\frac{k}{k+1}\int \eta\bk(\W)\dd\sig,
\]
where the last equality follows from Euler homogeneity. Adding the first term proves \eqref{first-var-M}.
\end{proof}

Let
\begin{equation*}
J_\varphi(s)
:=\int_1^s\frac{d\tau}{\tau\varphi(1/\tau)},
\qquad s>0.
\end{equation*}
Define the energy
\begin{equation}\label{energy}
\mathcal E_{f,\varphi}[h]
:=\mathcal M_k[h]
-\int_{\Ct}f\,J_\varphi(h/\ellc)\dd\sig.
\end{equation}

\begin{proposition}[Euler--Lagrange equation]\label{prop:EL}
A positive capillary function $h$ is a critical point of \eqref{energy} with respect to Robin variations if and only if it solves \eqref{main-eq}.
\end{proposition}

\begin{proof}
Because
\[
J_\varphi'(s)=\frac1{s\varphi(1/s)},
\]
the first variation of the second term in \eqref{energy} is
\[
\int_{\Ct}\frac{f}{h\varphi(\ellc/h)}\eta\dd\sig.
\]
Combining this with Lemma \ref{lem:first-var-M} gives
\[
D\mathcal E_{f,\varphi}[h](\eta)
=\int_{\Ct}\eta
\left[
\bk(\W[h])-\frac{f}{h\varphi(\ellc/h)}
\right]\dd\sig,
\]
which is equivalent to \eqref{main-eq}.
\end{proof}

At a critical point, every Robin variation can be written as $\eta=hw$ with $w_\mu=0$. If $S_k^{ij}=\partial\bk/\partial\W_{ij}$, a second variation calculation gives
\begin{equation}\label{second-var}
\begin{aligned}
D^2\mathcal E[h](hw,hw)
={}&-\int_{\Ct}h^2S_k^{ij}w_iw_j\dd\sig\\
&+\int_{\Ct}h\bk(\W)
\big(k+1-a_\varphi(\ellc/h)\big)w^2\dd\sig.
\end{aligned}
\end{equation}
Thus the supercritical condition \eqref{SC} is exactly the condition that makes the variational Hessian strictly negative in the admissible class. This observation will reappear as coercivity of the linearized operator in Section 5.

\subsection{Dilation, the critical weight, and the supercritical class}

The same quantity that appears in the variational Hessian is the infinitesimal dilation exponent. For $h_c=ch$ one has
\[
\W[h_c]=c\W[h],
\]
and therefore
\begin{equation*}
\mathcal D_{\varphi,k}[h_c]
:=h_c\varphi(\ellc/h_c)\bk(\W[h_c])
=c^{k+1}h\varphi(r/c)\bk(\W[h]),
\qquad r=\ellc/h.
\end{equation*}
Hence
\begin{equation*}
\left.\frac{d}{d\log c}\right|_{c=1}
\log \mathcal D_{\varphi,k}[h_c]
=k+1-a_\varphi(r)=-\Lambda_\varphi(r).
\end{equation*}
This gives a geometric explanation for the threshold $k+1$: in the supercritical regime the density is strictly decreasing under dilation, while the critical equation is dilation invariant.

\begin{proposition}[Critical and logarithmically convex weight classes]\label{prop:weight-class}
If $a_\varphi\equiv k+1$, then $\varphi(t)=Ct^{k+1}$ for some $C>0$. If $\varphi_0$ and $\varphi_1$ satisfy
\[
a_{\varphi_j}(t)\ge k+1+\delta\qquad(j=0,1),
\]
then the logarithmic interpolation
\[
\varphi_s(t)=\varphi_0(t)^{1-s}\varphi_1(t)^s
\]
satisfies the same inequality for every $s\in[0,1]$.
\end{proposition}

\begin{proof}
The first assertion follows by integrating the ODE $t(\log\varphi)'=k+1$. For the second,
\[
a_{\varphi_s}=(1-s)a_{\varphi_0}+sa_{\varphi_1}.
\]
\end{proof}

\section{A Priori Estimates}

Throughout this section $h$ is a positive smooth $k$-admissible solution of \eqref{main-eq}. We assume the supercritical condition \eqref{SC}. Unless otherwise stated, the constants in Propositions \ref{prop:C0}--\ref{prop:C2-boundary} may depend on $n,k,\theta,\delta$, positive lower and upper bounds for $f$, a fixed $C^2$ norm of $f$, and the values of $\varphi$ and its first two derivatives on the compact interval for $\ellc/h$ determined below, but not on the particular solution.

\subsection{\texorpdfstring{$C^0$}{C0} estimate}

Set
\[
u=\frac{h}{\ellc},
\qquad g=\frac{f}{\ellc}.
\]
By \eqref{u-Neumann}, $u_\mu=0$, while \eqref{Btheta} and \eqref{main-eq} give
\begin{equation}\label{u-equation}
u\varphi(1/u)\bk(\mathcal B_\theta[u])=g.
\end{equation}
Define
\begin{equation*}
\Gamma_{\varphi,k}(s):=s^{k+1}\varphi(1/s).
\end{equation*}
Then
\begin{equation*}
\frac{d\log\Gamma_{\varphi,k}}{d\log s}
=k+1-a_\varphi(1/s)\le-\delta.
\end{equation*}
Hence $\Gamma_{\varphi,k}$ is strictly decreasing.

\begin{proposition}[$C^0$ estimate]\label{prop:C0}
Let
\[
g_-:=\min_{\overline\Ct}g,
\qquad
g_+:=\max_{\overline\Ct}g.
\]
Then
\begin{equation*}
\Gamma_{\varphi,k}^{-1}(g_+)
\le\frac{h}{\ellc}
\le\Gamma_{\varphi,k}^{-1}(g_-).
\end{equation*}
In particular,
\[
0<c_0\le h\le C_0,
\qquad
0<c_0\le\frac{\ellc}{h}\le C_0.
\]
\end{proposition}

\begin{proof}
Let $x_M$ be a maximum point of $u$. Because $u_\mu=0$, the usual second derivative test applies also when $x_M\in\partial\Ct$; thus
\[
\nabla u(x_M)=0,
\qquad
\nabla^2u(x_M)\le0.
\]
By \eqref{Btheta},
\[
\mathcal B_\theta[u](x_M)\le u(x_M)\sig.
\]
Since $\bk$ is elliptic on $\Gamma_k$,
\[
\bk(\mathcal B_\theta[u](x_M))\le u(x_M)^k.
\]
Equation \eqref{u-equation} gives
\[
g(x_M)\le u(x_M)^{k+1}\varphi(1/u(x_M))
=\Gamma_{\varphi,k}(u_{\max}).
\]
Thus $g_-\le\Gamma_{\varphi,k}(u_{\max})$ and, because $\Gamma_{\varphi,k}$ is decreasing,
\[
u_{\max}\le\Gamma_{\varphi,k}^{-1}(g_-).
\]
The minimum point gives $\mathcal B_\theta[u]\ge u\sig$ and hence
\[
g_+\ge\Gamma_{\varphi,k}(u_{\min}),
\]
which yields the lower bound. Since $\ellc$ is bounded above and below by positive constants on an acute cap, the remaining estimates follow.
\end{proof}

The same argument gives spherical sub- and supersolutions. Indeed, for $h=c\ellc$,
\[
\W[c\ellc]=c\sig,
\qquad
h\varphi(\ellc/h)\bk(\W[h])
=\ellc\,\Gamma_{\varphi,k}(c).
\]
This observation is useful along the continuity path in Section 5.

\subsection{\texorpdfstring{$C^1$}{C1} estimate}

Let
\[
v=\log h.
\]
Then
\[
\W[h]=hB[v],
\qquad
B[v]:=\nabla^2v+\nabla v\otimes\nabla v+\sig.
\]
Equation \eqref{main-eq} can be written as
\begin{equation}\label{B-equation}
\bar\sigma_k(B[v])^{1/k}=\widetilde f,
\qquad
\widetilde f
:=\left[
\frac{f}{h^{k+1}\varphi(\ellc/h)}
\right]^{1/k}.
\end{equation}
Set $r=\ellc/h$. A direct differentiation gives the decisive identity
\begin{equation}\label{grad-ftilde}
\nabla\log\widetilde f
=\frac1k\nabla\log f
+\frac{a_\varphi(r)-k-1}{k}\nabla v
-\frac{a_\varphi(r)}k\nabla\log\ellc.
\end{equation}
The coefficient of $\nabla v$ satisfies
\begin{equation}\label{p0-positive}
\frac{a_\varphi(r)-k-1}{k}\ge\frac\delta k>0.
\end{equation}

\begin{proposition}[$C^1$ estimate]\label{prop:C1}
There exists $C_1$ such that
\begin{equation}\label{C1-estimate}
\|\nabla\log h\|_{L^\infty(\Ct)}\le C_1.
\end{equation}
Consequently $\|\nabla h\|_{L^\infty}\le C$.
\end{proposition}

\begin{proof}
We follow the logarithmic-gradient argument for the supercritical power equation, recording explicitly the only place where the Orlicz weight changes the calculation. Put
\[
\widetilde\ell:=\frac{\ellc}{\sin\theta\cos\theta},
\qquad
\omega:=v-\cot\theta\,\widetilde\ell,
\qquad
\Psi:=|\nabla\omega|^2.
\]
From \eqref{ell-identities},
\[
\widetilde\ell|_{\partial\Ct}=\tan\theta,
\qquad
\nabla\widetilde\ell|_{\partial\Ct}=\mu,
\qquad
\nabla^2\widetilde\ell
=\frac{1-\ellc}{\sin\theta\cos\theta}\,\sig.
\]
Since $v_\mu=\cot\theta$, we have
\begin{equation}\label{omega-Neumann}
\omega_\mu=0\qquad\text{on }\partial\Ct.
\end{equation}
Let
\[
\mathcal G(B):=\bar\sigma_k(B)^{1/k},
\qquad
\mathcal G^{ij}:=\frac{\partial\mathcal G}{\partial B_{ij}}.
\]
The operator is elliptic and concave on $\Gamma_k$, homogeneous of degree one, and satisfies
\begin{equation}\label{G-trace-lower}
\sum_i\mathcal G^{ii}\ge c_{n,k}>0.
\end{equation}

Let $x_0$ be a maximum point of $\Psi$. If $x_0\in\partial\Ct$, choose an orthonormal frame with $e_n=\mu$. Differentiating \eqref{omega-Neumann} tangentially and using the total umbilicity of $\partial\Ct\subset\Ct$ gives
\[
\omega_{\alpha n}=-\cot\theta\,\omega_\alpha.
\]
Because $e_n$ is the outward conormal and $\Psi$ has a maximum relative to $\overline\Ct$,
\[
0\le \Psi_n
=2\sum_{\alpha=1}^{n-1}\omega_\alpha\omega_{\alpha n}
=-2\cot\theta\sum_{\alpha=1}^{n-1}\omega_\alpha^2.
\]
Since $\cot\theta>0$, all tangential derivatives vanish, while \eqref{omega-Neumann} gives the normal derivative zero. Hence $|\nabla\omega|(x_0)=0$ and the estimate follows in the boundary case.

Assume now that $x_0$ is interior. Diagonalize $B[v]$ at $x_0$. If $|\nabla\omega|(x_0)$ is bounded there is nothing to prove, so assume it is larger than a fixed geometric constant. At the maximum,
\begin{equation}\label{Psi-first}
0=(\log\Psi)_i
=\frac{2\omega_s\omega_{si}}{\Psi}
\end{equation}
and
\begin{equation}\label{Psi-second}
0\ge \mathcal G^{ij}(\log\Psi)_{ij}
=\frac{2}{\Psi}\mathcal G^{ij}
\left(\omega_s\omega_{sij}+\omega_{si}\omega_{sj}\right).
\end{equation}
The spherical commutation formula is
\[
v_{sij}=v_{ijs}+v_s\delta_{ij}-v_j\delta_{si}.
\]
Differentiating \eqref{B-equation} gives
\[
\mathcal G^{ij}\big(v_{ijs}+2v_{si}v_j\big)=\widetilde f_s.
\]
Multiplying by $\omega_s$ and using \eqref{grad-ftilde},
\begin{align*}
\mathcal G^{ij}\omega_s\big(v_{ijs}+2v_{si}v_j\big)
&=\widetilde f_s\omega_s \notag\\
&=\widetilde f\left[
\frac1k(\log f)_s
+p_0(r)v_s
-\frac{a_\varphi(r)}k(\log\ellc)_s
\right]\omega_s,
\end{align*}
where $p_0(r)=(a_\varphi(r)-k-1)/k$. Since
\[
v_s=\omega_s+\cot\theta\,\widetilde\ell_s,
\qquad
p_0(r)\ge\delta/k,
\]
and Proposition \ref{prop:C0} confines $r$ to a compact subinterval of $(0,\infty)$, all coefficients other than the quadratic term are uniformly bounded. Moreover $\widetilde f$ is bounded above and below by positive constants. Therefore
\begin{equation}\label{C1-RHS-lower}
\widetilde f_s\omega_s
\ge c\Psi-C(1+\sqrt\Psi).
\end{equation}
This is the precise Orlicz substitute for the constant supercritical coefficient in the power-law proof.

It remains only the operator algebra, which is unchanged. Write
\[
\omega_{si}=B_{si}-v_sv_i-\delta_{si}-\cot\theta\,\widetilde\ell_{si}.
\]
From \eqref{Psi-first}, $\omega_s\omega_{si}=0$. At the diagonal point this gives
\[
\omega_iB_{ii}
=\cot\theta\,\widetilde\ell_{si}\omega_s
+\omega_i+(\omega_sv_s)v_i.
\]
Using this identity to complete squares, together with
\[
\widetilde\ell_{ij}
=\frac{1-\ellc}{\sin\theta\cos\theta}\delta_{ij},
\qquad
\mathcal G^{ij}B_{ij}=\widetilde f,
\]
yields
\begin{align}
&2\sum_{s,i}\mathcal G^{ii}\omega_{si}^2
-4\sum_{s,i}\mathcal G^{ii}\omega_sv_{si}v_i \notag\\
&\qquad\ge
4\sum_i\mathcal G^{ii}v_i^2
-C(1+|\nabla v|)\left(1+\sum_i\mathcal G^{ii}\right).
\label{C1-square-est}
\end{align}
Substituting the commutation identity, \eqref{C1-RHS-lower}, and \eqref{C1-square-est} into \eqref{Psi-second}, and using that $|\nabla v|$ and $|\nabla\omega|$ differ by a fixed bounded vector, gives for $|\nabla\omega|$ sufficiently large
\[
0\ge
2\sum_i\mathcal G^{ii}
-\frac{C}{|\nabla\omega|}
\left(1+\sum_i\mathcal G^{ii}\right).
\]
By \eqref{G-trace-lower}, this is impossible when $|\nabla\omega|$ is large. Thus $|\nabla\omega|\le C$ on $\overline\Ct$, and the fixed identity $v=\omega+\cot\theta\widetilde\ell$ proves \eqref{C1-estimate}.
\end{proof}

\begin{remark}
The key point is not the exact choice of the boundary correction $\omega$, but the strictly positive coefficient \eqref{p0-positive}. It is the Orlicz replacement of the supercritical power coefficient $p-k-1$.
\end{remark}

\subsection{Interior \texorpdfstring{$C^2$}{C2} estimate}

Set
\begin{equation*}
F(\W):=\bk(\W)^{1/k},
\qquad
R:=\left[\frac{f}{h\varphi(\ellc/h)}\right]^{1/k},
\end{equation*}
so that
\begin{equation*}
F(\W)=R.
\end{equation*}
Let
\[
P:=\max_{\xi\in\overline\Ct,\ |\eta|=1}\W(\eta,\eta).
\]
Suppose first that the maximum is achieved at an interior point. Choose a local orthonormal frame so that $\W$ is diagonal and $P=\W_{11}$ there. The spherical commutation formula, concavity of $F$, and $\nabla\W_{11}=0$ give the standard inequality
\begin{equation}\label{interior-Wmax}
0\ge R_{11}+P\sum_iF^{ii}-R.
\end{equation}
This is the same principal-radius maximum principle used in the power problem \cite{MWWLp}.

Let
\[
L=\log R,
\qquad q=\log\ellc.
\]
Since $v=\log h$ and $r=\ellc/h$, direct differentiation yields
\begin{equation*}
\begin{aligned}
L_{11}=\frac1k\Big[
&(\log f)_{11}
+\frac{a_\varphi(r)-1}{h}P
+1-\frac{a_\varphi(r)}{\ellc}\\
&-(a_\varphi(r)-1)v_1^2
+a_\varphi(r)q_1^2
-r a_\varphi'(r)(q_1-v_1)^2
\Big].
\end{aligned}
\end{equation*}
By Propositions \ref{prop:C0}--\ref{prop:C1}, all terms except the one containing $P$ are bounded. Moreover,
\[
a_\varphi(r)-1\ge k+\delta.
\]
Hence
\[
L_{11}\ge cP-C.
\]
Since $R$ is bounded above and below by positive constants,
\begin{equation*}
R_{11}=R(L_{11}+L_1^2)\ge cP-C.
\end{equation*}
Substituting into \eqref{interior-Wmax} gives $P\le C$ at an interior maximum.

\begin{proposition}[Interior second-derivative estimate]\label{prop:C2-interior}
If the maximum eigenvalue of $\W[h]$ is attained in the interior, then
\[
\lambda_{\max}(\W[h])\le C.
\]
\end{proposition}

\subsection{Boundary \texorpdfstring{$C^2$}{C2} estimate}

The remaining task is to control the largest eigenvalue when it is attained on $\partial\Ct$. The tangential second derivatives reduce to the double-normal derivative through the capillary boundary identities. More precisely, the global reduction in the acute cap has the form
\begin{equation}\label{C2-reduction}
\sup_{\Ct}|\nabla^2h|
\le C\left(1+M\right),
\qquad
M:=\sup_{\partial\Ct}|h_{\mu\mu}|.
\end{equation}
This is the same reduction used in the supercritical capillary $L_p$ Hessian estimate; see \cite[Lemma~3.6(1), especially Case~2 of its proof]{MWWLp}. Its proof uses the vanishing of the mixed components in \eqref{W-mixed-zero}, the boundary identity \eqref{W-normal-derivative}, and the sign $\cot\theta>0$, so it applies unchanged to the present right-hand side once the $C^0$--$C^1$ bounds are available.

To estimate $M$, let $d$ be the distance to $\partial\Ct$ in a sufficiently small collar and set
\begin{equation*}
\zeta=e^{-d}-1.
\end{equation*}
For $0<\theta<\pi/2$, the boundary of the spherical cap is strictly convex and one has, after shrinking the collar,
\begin{equation*}
\nabla^2\zeta\ge c_\theta\sig
\end{equation*}
for some $c_\theta>0$. This is the source of the acute-angle restriction.

The differentiated right-hand side satisfies
\begin{equation}\label{grad-R}
\nabla\log R
=\frac1k\left[
\nabla\log f
+(a_\varphi(r)-1)\nabla\log h
-a_\varphi(r)\nabla\log\ellc
\right].
\end{equation}
Thus Propositions \ref{prop:C0}--\ref{prop:C1} imply
\begin{equation}\label{R-C1}
0<c\le R\le C,
\qquad
|\nabla R|\le C.
\end{equation}

\begin{proposition}[Boundary double-normal estimate]\label{prop:C2-boundary}
There is a constant $C$ such that
\begin{equation}\label{M-bound}
\sup_{\partial\Ct}|h_{\mu\mu}|\le C.
\end{equation}
Consequently,
\begin{equation}\label{global-C2}
\|h\|_{C^2(\overline\Ct)}\le C.
\end{equation}
\end{proposition}

\begin{proof}
We give the barrier computation in the form needed for $F=\bar\sigma_k^{1/k}$. Let
\[
\Omega_{\delta_0}:=\{\xi\in\Ct:0\le d(\xi)\le\delta_0\},
\]
where $d$ is the inward distance from $\partial\Ct$. In this collar,
\[
\zeta=e^{-d}-1,
\qquad
\zeta|_{\partial\Ct}=0,
\qquad
\nabla\zeta|_{\partial\Ct}=\mu,
\]
and
\[
\zeta_{ij}
=e^{-d}\left[
\cot(\theta-d)(\sig_{ij}-d_id_j)+d_id_j
\right].
\]
After decreasing $\delta_0$ there is $\eps_0>0$, depending only on $n,k,\theta$, such that
\begin{equation*}
F^{ij}\zeta_{ij}
\ge \eps_0\left(\sum_iF^{ii}+1\right)
\qquad\text{in }\Omega_{\delta_0}.
\end{equation*}
All second and third derivatives of $d$ are uniformly bounded in this fixed collar.

Set
\[
\mathcal N[h]:=\langle\nabla h,\nabla d\rangle+\cot\theta\,h.
\]
The Robin condition gives $\mathcal N[h]=0$ on $\partial\Ct$. For the upper estimate define
\begin{equation}\label{upper-barrier}
\overline Q
:=-(1+\beta d)\mathcal N[h]
-\left(B_1+\frac M2\right)\zeta,
\end{equation}
where $\beta>1$ and $B_1>1$ are constants to be chosen. Clearly $\overline Q=0$ on $\partial\Ct$.

Suppose $\overline Q$ has a negative minimum at an interior point $x_*$. Choose an orthonormal frame at $x_*$ in which $\W$ and $F^{ij}$ are diagonal. Differentiating \eqref{upper-barrier}, using the spherical commutation identity
\[
h_{kij}=\W_{ijk}-h_j\delta_{ki},
\]
and contracting with $F^{ij}$ gives
\begin{align}
F^{ij}\overline Q_{ij}
\le{}&-2\beta\sum_iF^{ii}\W_{ii}d_i^2
-2(1+\beta d)\sum_iF^{ii}\W_{ii}d_{ii}\notag\\
&-\left[
\left(B_1+\frac M2\right)\eps_0-C(\beta+1)
\right]
\left(\sum_iF^{ii}+1\right).
\label{barrier-core}
\end{align}
The only differentiated-equation term is
\[
F^{ij}\W_{ijs}d_s=R_sd_s,
\]
which is uniformly bounded by \eqref{R-C1}. Thus every constant in \eqref{barrier-core} is already controlled by the $C^0$--$C^1$ estimates and the fixed geometry.

Put
\[
\mathcal B:=\left\{i:\beta d_i^2<\frac1n\right\},
\qquad
\mathcal G:=\left\{i:\beta d_i^2\ge\frac1n\right\}.
\]
Since $|\nabla d|=1$, $\mathcal G$ is nonempty; relabel so that $1\in\mathcal G$. At the minimum, $\overline Q_i=0$. Using
\[
(\mathcal N[h])_i
=(\W_{ii}-h)d_i+h_s d_{si}+\cot\theta\,h_i
\]
in the diagonal frame, the identity $\overline Q_i=0$ becomes
\[
(1+\beta d)\big[(\W_{ii}-h)d_i+h_s d_{si}+\cot\theta\,h_i\big]
=\left(B_1+\frac M2\right)e^{-d}d_i-\beta d_i\mathcal N[h].
\]
For $i\in\mathcal G$ we have $|d_i|\ge(n\beta)^{-1/2}$. The $C^0$--$C^1$ bounds and the fixed collar geometry therefore give
\[
\left|
\W_{ii}-
\frac{(B_1+M/2)e^{-d}}{1+\beta d}
\right|
\le C(1+\beta),
\qquad i\in\mathcal G.
\]
Taking
\[
\delta_0\le\min\left\{\frac\theta2,\frac1{9\beta+10}\right\}
\]
and then choosing $B_1\ge C_1\beta$ gives
\begin{equation}\label{Wgood-bounds}
\frac{3B_1}{5}+\frac{2M}{5}
\le \W_{ii}
\le \frac{6B_1}{5}+\frac M2,
\qquad i\in\mathcal G.
\end{equation}
In particular $\W_{11}>0$. Using the concavity of $d$, the definition of $\mathcal B$, and the homogeneity identity $\sum_iF^{ii}\W_{ii}=R$, inequality \eqref{barrier-core} reduces to
\begin{align}
F^{ij}\overline Q_{ij}
\le{}&-\frac{2\beta}{n}F^{11}\W_{11}
+C_0\sum_{\W_{ii}>0}F^{ii}\W_{ii}\notag\\
&-\left[
\left(B_1+\frac M2\right)\eps_0-C_0(\beta+1)
\right]
\left(\sum_iF^{ii}+1\right).
\label{barrier-reduced}
\end{align}
Order the remaining entries so that $\W_{22}\ge\cdots\ge\W_{nn}$. We consider three cases.

If $\W_{nn}\ge0$, all eigenvalues are nonnegative and
\[
\sum_{\W_{ii}>0}F^{ii}\W_{ii}=R\le C.
\]
Thus, for fixed $\beta$, a sufficiently large $B_1$ makes the right-hand side of \eqref{barrier-reduced} negative.

Suppose next that $\W_{nn}<0$ but $-\W_{nn}\le\eta\W_{11}$, where
\[
\eta:=\frac{\eps_0}{10(C_0+1)}.
\]
Since every negative eigenvalue is at least $\W_{nn}$,
\[
\sum_{\W_{ii}>0}F^{ii}\W_{ii}
=R-\sum_{\W_{ii}<0}F^{ii}\W_{ii}
\le C+\eta\W_{11}\sum_iF^{ii}.
\]
Together with the upper bound in \eqref{Wgood-bounds}, this term is absorbed by a fixed fraction of
\[
-\eps_0\left(B_1+\frac M2\right)\sum_iF^{ii},
\]
and again $B_1$ sufficiently large gives $F^{ij}\overline Q_{ij}<0$.

Finally suppose $\W_{nn}<0$ and $-\W_{nn}>\eta\W_{11}$. The reduction \eqref{C2-reduction} and \eqref{Wgood-bounds} imply
\[
\W_{22}\le C_4(1+M),
\qquad
\W_{11}\ge c_4(1+M).
\]
For $F=\bar\sigma_k^{1/k}$, the G\aa rding-cone derivative-separation estimate used in the proof of \cite[Lemma~3.7]{MWWLp} gives
\[
F^{11}\ge\eps_1\sum_iF^{ii}
\]
for some $\eps_1>0$ depending only on the preceding ratio constants; when $k=1$ all $F^{ii}$ are equal. Hence the first term of \eqref{barrier-reduced} is at most
\[
-c\beta(1+M)\sum_iF^{ii},
\]
while the positive-eigenvalue term is at most $C(1+M)\sum_iF^{ii}$. Choose $\beta$ sufficiently large and then $B_1$ sufficiently large. This again gives $F^{ij}\overline Q_{ij}<0$.

All three cases contradict the existence of an interior negative minimum. With $\beta$ and $\delta_0$ fixed, increasing $B_1$ once more makes $\overline Q>0$ on the artificial boundary $\{d=\delta_0\}$. Therefore
\begin{equation}\label{Qbar-positive}
\overline Q\ge0\qquad\text{in }\Omega_{\delta_0}.
\end{equation}
Let $x^+\in\partial\Ct$ maximize $h_{\mu\mu}$. Since $\overline Q=0$ on the true boundary and \eqref{Qbar-positive} holds on the interior side, the outward derivative is nonpositive. Using $\zeta_\mu=1$ and the boundary identities gives
\[
0\ge\overline Q_\mu(x^+)
\ge h_{\mu\mu}(x^+)-\left(B_1+\frac M2\right)-C,
\]
so
\begin{equation}\label{hmunu-upper}
\sup_{\partial\Ct}h_{\mu\mu}\le C+\frac M2.
\end{equation}

For the lower estimate define
\begin{equation*}
\underline Q
:=-(1+\underline\beta d)\mathcal N[h]
+\left(B_2+\frac M2\right)\zeta.
\end{equation*}
The same calculation with reversed inequalities, enlarging the fixed constants if necessary, excludes an interior positive maximum and makes $\underline Q<0$ on the artificial boundary. Hence $\underline Q\le0$ in the collar. At a point $x^-\in\partial\Ct$ minimizing $h_{\mu\mu}$,
\begin{equation}\label{hmunu-lower}
\inf_{\partial\Ct}h_{\mu\mu}\ge-C-\frac M2.
\end{equation}
By the definition of $M$, \eqref{hmunu-upper}--\eqref{hmunu-lower} imply
\[
M\le C+\frac M2,
\]
so $M\le2C$. This proves \eqref{M-bound}; \eqref{global-C2} follows from \eqref{C2-reduction}.
\end{proof}

\begin{remark}
The boundary $C^2$ estimate is the only place in the admissible existence theory where the acute-angle sign $\cot\theta>0$ is essential. No extra Orlicz structural condition beyond \eqref{SC} is used: the Orlicz dependence enters only through the bounded quantities in \eqref{grad-R}.
\end{remark}

\subsection{Higher regularity and uniform ellipticity}

The global $C^2$ estimate bounds every eigenvalue of $\W[h]$ from above in absolute value. On the other hand,
\[
\bk(\W[h])=\frac{f}{h\varphi(\ellc/h)}
\]
has a uniform positive lower bound by the $C^0$ estimate and the positivity of $f$ and $\varphi$ on the compact range of $\ellc/h$. Lemma \ref{lem:Gamma-compact} therefore gives a compact set
\[
K\Subset\Gamma_k
\]
such that
\[
\lambda(\W[h](x))\in K\qquad\text{for every }x\in\overline\Ct.
\]
Consequently the coefficients $F^{ij}(\W[h])$ are uniformly positive definite and uniformly bounded. Since $F=\bar\sigma_k^{1/k}$ is concave on $\Gamma_k$, the interior Evans--Krylov estimate applies. At the boundary the condition is the fixed linear oblique condition
\[
h_\mu-\cot\theta h=0,
\]
so the standard boundary $C^{2,\alpha}$ estimate for uniformly elliptic concave equations with oblique boundary data gives
\begin{equation*}
\|h\|_{C^{2,\alpha}(\overline\Ct)}\le C_\alpha.
\end{equation*}
After this estimate is available, differentiating the equation and applying Schauder estimates inductively yields smoothness. For each $m\ge1$, however, the $C^{m+2,\alpha}$ estimate also depends on the corresponding $C^{m,\alpha}$ norm of $f$ and on finitely many derivatives of $\varphi$ on the compact range of $\ellc/h$. Thus higher-order estimates are uniform only when these higher-order data norms are uniform.

\section{Proof of the Main Theorem by the Continuity Method}

\subsection{The continuity path}

Fix the Orlicz weight $\varphi$ and define
\begin{equation}\label{data-path}
f_t:=(1-t)\varphi(1)\ellc+t f,
\qquad 0\le t\le1.
\end{equation}
Every $f_t$ is positive. At $t=0$ the spherical cap support function
\[
h_0=\ellc
\]
is an exact solution because $\W[\ellc]=\sig$ and $\bk(\sig)=1$.

Let
\[
\mathcal I:=\left\{t\in[0,1]:
\text{\eqref{main-eq} with datum $f_t$ admits a positive smooth $k$-admissible solution}
\right\}.
\]
Then $0\in\mathcal I$. We prove that $\mathcal I$ is open and closed.

\subsection{Linearization, self-adjointness, and openness}

Write the equation in logarithmic form
\begin{equation}\label{G-map}
\mathcal G_t(h)
:=\log\bk(\W[h])
+\log h
+\log\varphi(\ellc/h)
-\log f_t=0.
\end{equation}
Let
\[
A^{ij}:=\frac{\partial\log\bk}{\partial\W_{ij}}
=\frac{S_k^{ij}}{\bk(\W)}.
\]
If $\eta$ is a Robin variation, write
\[
\eta=hw.
\]
Since both $h$ and $\eta$ satisfy the Robin condition,
\begin{equation*}
w_\mu=0.
\end{equation*}
Using Euler's identity $A^{ij}\W_{ij}=k$, the linearization of \eqref{G-map} is
\begin{equation*}
\mathscr L_h w
=hA^{ij}w_{ij}
+2A^{ij}h_iw_j
+\big(k+1-a_\varphi(\ellc/h)\big)w.
\end{equation*}
Thus the supercritical condition supplies a strictly negative zeroth-order term:
\begin{equation}\label{zeroth-negative}
k+1-a_\varphi(\ellc/h)\le-\delta.
\end{equation}

There is also a useful self-adjoint formulation. Define
\begin{equation*}
\dd\mu_h:=h\bk(\W[h])\dd\sig.
\end{equation*}
Since $S_k^{ij}$ is divergence-free, $\W_{\alpha n}=0$ at the boundary, and $w_\mu=0$, integration by parts gives for Neumann functions $u,w$
\begin{equation}\label{self-adjoint}
\begin{aligned}
\int_{\Ct}u\mathscr L_hw\dd\mu_h
={}&-\int_{\Ct}h^2S_k^{ij}u_iw_j\dd\sig\\
&+\int_{\Ct}
\big(k+1-a_\varphi(\ellc/h)\big)uw\dd\mu_h.
\end{aligned}
\end{equation}
In particular,
\begin{equation}\label{coercivity}
-\int w\mathscr L_hw\dd\mu_h
\ge
\int h^2S_k^{ij}w_iw_j\dd\sig
+\delta\int w^2\dd\mu_h.
\end{equation}
This is precisely the second-variation coercivity seen in \eqref{second-var}.

\begin{proposition}[Self-adjoint spectral gap]\label{prop:spectral-gap}
The symmetric quadratic form determined by \eqref{self-adjoint} on $H^1(\Ct)$ is closed and coercive. Its Friedrichs realization defines a self-adjoint operator $-\mathscr L_h$ on $L^2(\Ct,d\mu_h)$ with the natural Neumann boundary condition, and
\begin{equation}\label{spectral-gap}
\lambda_1(-\mathscr L_h)\ge\delta.
\end{equation}
Equivalently, every smooth Neumann function $w$ satisfies
\[
\|w\|_{L^2(d\mu_h)}^2
\le\frac1\delta
\left(-\int_{\Ct}w\mathscr L_hw\,d\mu_h\right).
\]
\end{proposition}

\begin{proof}
Symmetry follows from \eqref{self-adjoint}. The fixed smooth admissible solution makes $h^2S_k^{ij}$ uniformly elliptic and the measures $d\mu_h$ and $d\sig$ uniformly equivalent. Together with \eqref{SC}, \eqref{coercivity} therefore gives a closed coercive form on $H^1(\Ct)$. The Friedrichs representation theorem yields the self-adjoint realization, and the Rayleigh quotient gives \eqref{spectral-gap}.
\end{proof}

\begin{proposition}[Linearized isomorphism]\label{prop:isomorphism}
At every positive smooth $k$-admissible solution, the linearized Robin operator is an isomorphism between the standard H\"older spaces with the boundary condition built into the domain. Consequently $\mathcal I$ is open.
\end{proposition}

\begin{proof}
For Neumann functions $u,w$ define
\[
\mathfrak a(u,w)
:=\int_{\Ct}h^2S_k^{ij}u_iw_j\dd\sig
+\int_{\Ct}\big(a_\varphi(\ellc/h)-k-1\big)uw\dd\mu_h.
\]
By \eqref{self-adjoint}, this is the weak bilinear form of $-\mathscr L_h$. Since the fixed solution is $k$-admissible and smooth on the compact cap, its eigenvalues lie in a compact subset of $\Gamma_k$. Thus $S_k^{ij}$ is uniformly positive definite. The functions $h$, $\bk(\W)$ and their reciprocals are bounded, and \eqref{SC} gives
\[
\mathfrak a(w,w)
\ge c\|\nabla w\|_{L^2(\Ct)}^2
+\delta\int_{\Ct}w^2\dd\mu_h
\ge c_0\|w\|_{H^1(\Ct)}^2.
\]
Hence Lax--Milgram gives a unique weak solution $w\in H^1(\Ct)$ of
\[
-\mathscr L_h w=g
\]
for every $g\in H^{-1}(\Ct)$, with the natural conormal boundary condition associated with the form. If $g\in C^{m,\alpha}$, elliptic boundary regularity makes $w$ classical. In a boundary frame with $e_n=\mu$, \eqref{W-mixed-zero} makes $\W$ block diagonal and hence $S_k^{\alpha n}=0$. The conormal condition is therefore
\[
h^2S_k^{nn}w_n=0.
\]
Since $S_k^{nn}>0$ on $\Gamma_k$, it is exactly $w_\mu=0$.

Uniform ellipticity and this linear Neumann boundary condition give the Schauder estimate
\[
\|w\|_{C^{m+2,\alpha}}
\le C\big(\|g\|_{C^{m,\alpha}}+\|w\|_{C^0}\big).
\]
The spectral gap, or equivalently the maximum principle using \eqref{zeroth-negative}, controls the $C^0$ term and shows that the homogeneous problem has only the zero solution. Therefore
\[
\mathscr L_h:
\{w\in C^{m+2,\alpha}:w_\mu=0\}
\longrightarrow C^{m,\alpha}
\]
is bijective with bounded inverse. Multiplication by the positive function $h$ identifies this Neumann space with the corresponding Robin variation space. The implicit function theorem then proves that $\mathcal I$ is open.
\end{proof}

\subsection{Closedness}

Let $t_j\in\mathcal I$ and $t_j\to t_\infty$, and let $h_j$ be the corresponding solutions. The path \eqref{data-path} has uniform positive lower and upper bounds and uniform $C^2$ norms. Hence Propositions \ref{prop:C0}--\ref{prop:C2-boundary} and the $C^{2,\alpha}$ estimate of Section 4.5 are uniform in $j$; in particular,
\[
\|h_j\|_{C^{2,\alpha}(\overline\Ct)}\le C.
\]
The global $C^2$ bound gives a common upper bound for $|\lambda(\W[h_j])|$, while the equation and the uniform $C^0$ estimate give
\[
\bk(\W[h_j])
=\frac{f_{t_j}}{h_j\varphi(\ellc/h_j)}
\ge c>0.
\]
Lemma \ref{lem:Gamma-compact} therefore supplies one compact set $K\Subset\Gamma_k$ such that
\[
\lambda(\W[h_j](x))\in K
\qquad\text{for every }j\text{ and }x\in\overline\Ct.
\]
Thus the equations are uniformly elliptic with constants independent of $j$.

After passing to a subsequence,
\[
h_j\longrightarrow h_\infty
\quad\text{in }C^{2,\beta}(\overline\Ct),
\qquad 0<\beta<\alpha.
\]
The equation and the Robin boundary condition pass to the limit. Since $K$ is compactly contained in $\Gamma_k$, the limiting tensor still satisfies
\[
\lambda(\W[h_\infty])\in K\Subset\Gamma_k,
\]
so admissibility cannot be lost. Since $f_{t_\infty}$ and $\varphi$ are smooth and the limiting equation is uniformly elliptic, standard Schauder bootstrapping applied directly to $h_\infty$ yields $h_\infty\in C^\infty$. Hence $t_\infty\in\mathcal I$, and $\mathcal I$ is closed.

Since $[0,1]$ is connected and $0\in\mathcal I$,
\[
\mathcal I=[0,1].
\]
This proves the existence of a positive smooth $k$-admissible solution of \eqref{main-eq}.

\subsection{Uniqueness and quantitative comparison}

The same supercritical slope that gives openness also gives a direct comparison principle.

\begin{proposition}[Comparison principle]\label{prop:comparison}
Let $h_1,h_2$ be positive $k$-admissible solutions with data $f_1,f_2$ and the same Orlicz weight $\varphi$. Then
\begin{equation}\label{comparison-stability}
\left\|\log\frac{h_1}{h_2}\right\|_{L^\infty(\Ct)}
\le\frac1\delta
\left\|\log\frac{f_1}{f_2}\right\|_{L^\infty(\Ct)}.
\end{equation}
In particular, solutions with the same datum coincide.
\end{proposition}

\begin{proof}
Let
\[
M=\max_{\overline\Ct}\frac{h_1}{h_2}.
\]
At a maximum point, the quotient has zero normal derivative if the point lies on the boundary, because both functions satisfy the same Robin condition. Hence
\[
\W[h_1]\le M\W[h_2].
\]
By ellipticity and homogeneity,
\[
\bk(\W[h_1])\le M^k\bk(\W[h_2]).
\]
Writing $r_2=\ellc/h_2$ and taking the ratio of the equations at the maximum point gives
\[
\frac{f_1}{f_2}
\le M^{k+1}\frac{\varphi(r_2/M)}{\varphi(r_2)}.
\]
If $M\ge1$, integrating \eqref{SC} between $r_2/M$ and $r_2$ yields
\[
\frac{\varphi(r_2/M)}{\varphi(r_2)}
\le M^{-(k+1+\delta)},
\]
so
\[
M^\delta\le\frac{f_2}{f_1}
\]
at that point. Thus
\[
M\le\exp\left(\frac1\delta
\left\|\log\frac{f_1}{f_2}\right\|_\infty\right).
\]
Applying the same argument to $h_2/h_1$ gives \eqref{comparison-stability}. If $f_1=f_2$, both maxima are at most one, so $h_1=h_2$.
\end{proof}

\subsection{Differentiation of the solution map and joint stability}

The coercive linearization also gives quantitative dependence on both the datum and the Orlicz weight. Let $s\mapsto(f_s,\varphi_s,h_s)$ be a smooth family of supercritical data and corresponding admissible solutions with a common margin $\delta>0$. Set
\[
w_s:=\partial_s\log h_s,
\qquad
g_s:=\partial_s\log f_s,
\qquad
\psi_s(t):=\partial_s\log\varphi_s(t),
\qquad
r_s:=\frac{\ellc}{h_s}.
\]
Differentiating the logarithmic equation gives
\begin{equation}\label{solution-derivative}
\mathscr L_{h_s}w_s
=g_s-\psi_s(r_s),
\qquad
(w_s)_\mu=0.
\end{equation}
The negative zeroth-order coefficient and the maximum principle imply
\begin{equation}\label{solution-Linf}
\|w_s\|_{L^\infty(\Ct)}
\le\frac1\delta
\|g_s-\psi_s(r_s)\|_{L^\infty(\Ct)}.
\end{equation}
Schauder theory applied to \eqref{solution-derivative} yields corresponding $C^{2,\alpha}$ estimates and smooth dependence on the data.

\begin{corollary}[Joint logarithmic stability]\label{cor:joint-stability}
Let $(f_0,\varphi_0)$ and $(f_1,\varphi_1)$ be two positive data pairs satisfying the same supercritical margin $\delta$, and suppose the logarithmic interpolation
\[
f_s=f_0^{1-s}f_1^s,
\qquad
\varphi_s=\varphi_0^{1-s}\varphi_1^s
\]
stays in a class for which the corresponding $r_s=\ellc/h_s$ lie in a common compact interval $I\Subset(0,\infty)$. Then
\begin{equation}\label{joint-stability}
\left\|\log\frac{h_1}{h_0}\right\|_{L^\infty}
\le\frac1\delta\left[
\left\|\log\frac{f_1}{f_0}\right\|_{L^\infty(\Ct)}
+
\left\|\log\frac{\varphi_1}{\varphi_0}\right\|_{L^\infty(I)}
\right].
\end{equation}
\end{corollary}

\begin{proof}
For the logarithmic interpolation, $g_s=\log(f_1/f_0)$ and $\psi_s=\log(\varphi_1/\varphi_0)$ are independent of $s$. Integrating \eqref{solution-Linf} from $s=0$ to $s=1$ proves \eqref{joint-stability}.
\end{proof}

The sign in \eqref{solution-derivative} also gives monotonicity: increasing $f$ with fixed $\varphi$ decreases $h$, while increasing $\varphi$ with fixed $f$ increases $h$.

\begin{proposition}[Open stability of strict convexity]\label{prop:strict-open}
Suppose $(f_0,\varphi_0)$ is supercritical with margin $\delta>0$ and its admissible solution satisfies
\[
\W[h_0]\ge\kappa\sig
\]
for some $\kappa>0$. Fix $0<\alpha<1$ and a compact interval $I\Subset(0,\infty)$ whose interior contains the range of $\ellc/h_0$. Then all sufficiently small positive smooth perturbations of $(f_0,\varphi_0)$ in $C^\alpha(\overline\Ct)\times C^2(I)$, preserving the same supercritical margin, have unique admissible solutions with $\ellc/h\in I$ and satisfying
\[
\W[h]\ge\frac\kappa2\sig.
\]
In particular, Theorem \ref{thm:mainC} holds.
\end{proposition}

\begin{proof}
Proposition \ref{prop:isomorphism} and the implicit function theorem give $C^{2,\alpha}$-continuous, indeed smooth, dependence of $h$ on $(f,\varphi)$ on the relevant compact $r$-interval. Since $h\mapsto\nabla^2h+h\sig$ is continuous from $C^{2,\alpha}$ to $C^{\alpha}$, the lower eigenvalue bound persists under sufficiently small perturbations.
\end{proof}

Thus, under the supercritical condition \eqref{SC}, the Robin problem \eqref{main-eq} has a unique positive smooth $k$-admissible solution.

\subsection{Full rank under the capillary Orlicz-convexity condition}

It remains to prove strict convexity under the additional assumptions of Theorem \ref{thm:mainA}. The key input is a non-homogeneous constant-rank identity adapted to the Orlicz right-hand side.

Write
\[
r=\frac{\ellc}{h},
\qquad
R=\left[\frac{f}{h\varphi(r)}\right]^{1/k},
\qquad
L=\log R,
\]
so that $F(\W)=\bk(\W)^{1/k}=R$. Define
\[
a(t)=a_\varphi(t),
\qquad
b(t)=a(t)-1,
\qquad
d(t)=t a'(t),
\]
\begin{equation*}
\mathcal A_\varphi(t):=b(t)+d(t)+\frac{b(t)^2}{k},
\end{equation*}
and
\[
G:=\log\frac{f}{\ellc},
\qquad q:=\log\ellc.
\]
Set
\begin{equation*}
\begin{aligned}
\mathscr S_{\varphi,f}(t)
={}&-\nabla^2G+\frac1k\dd G\otimes\dd G
+\left(k+\frac{b(t)}{\ellc}\right)\sig\\
&-\frac{b(t)^2}{\mathcal A_\varphi(t)}
\left(\dd q+\frac1k\dd G\right)^{\otimes2}.
\end{aligned}
\end{equation*}

\begin{lemma}[Orlicz square completion]\label{lem:square-completion}
Let $e$ be a unit eigenvector of $\W$ and put $\rho=\log(h/\ellc)$. Then
\begin{equation}\label{square-completion}
\begin{aligned}
k(L_e^2-L_{ee}+1)
={}&\mathcal A_\varphi(r)
\left[
\rho_e+\frac{b(r)}{\mathcal A_\varphi(r)}
\left(q_e+\frac1kG_e\right)
\right]^2\\
&+\mathscr S_{\varphi,f}(r)(e,e)
-\frac{b(r)}h\W(e,e).
\end{aligned}
\end{equation}
\end{lemma}

\begin{proof}
Since $r=e^{-\rho}$,
\[
kL=G-\rho-\log\varphi(e^{-\rho}).
\]
Hence
\[
kL_e=G_e+b\rho_e,
\qquad
kL_{ee}=G_{ee}+b\rho_{ee}-d\rho_e^2.
\]
The identities $\W(e,e)=h(v_{ee}+v_e^2+1)$ and \eqref{ell-identities} imply
\[
\rho_{ee}=\frac{\W(e,e)}h-\frac1{\ellc}-2q_e\rho_e-\rho_e^2.
\]
Substitution gives a quadratic polynomial in $\rho_e$ with leading coefficient $\mathcal A_\varphi$. Completing the square yields \eqref{square-completion}.
\end{proof}

If $\mathcal A_\varphi>0$ and $\mathscr S_{\varphi,f}\ge0$, \eqref{square-completion} gives
\[
L_e^2-L_{ee}+1\ge-C\W(e,e).
\]
To verify the non-homogeneous constant-rank hypothesis precisely, set
\[
\mathscr F(A,x):=F(A)-R(x),
\qquad F(A)=\bar\sigma_k(A)^{1/k}.
\]
For a unit eigenvector $e$ of $\W$, the curvature-adjusted inverse-concavity calculation in the proof of \cite[Theorem~1.2]{BISCRT} gives on the unit sphere
\[
\omega_{\mathscr F}(\W)(\nabla_e\W,e)
\ge 2R^{-1}R_e^2-R_{ee}+R-C\W(e,e).
\]
Since $R=e^L$,
\[
2R^{-1}R_e^2-R_{ee}+R
=R(L_e^2-L_{ee}+1),
\]
and the preceding square-completion inequality therefore yields
\[
\omega_{\mathscr F}(\W)(\nabla_e\W,e)
\ge-C\W(e,e).
\]
Thus the hypothesis of \cite[Theorem~1.6]{BISCRT} is satisfied (in fact in a stronger form than required). Because $\W$ is Codazzi and $F=\bar\sigma_k^{1/k}$ is strictly elliptic and inverse-concave on $\Gamma_k$, the Bryan--Ivaki--Scheuer theorem applies on the interior. Hence, if $\W\ge0$ and an interior rank degeneration occurs, $\W$ has constant rank. To justify \eqref{capillary-Minkowski-identity} in the semidefinite case, apply it first to $h_\varepsilon=h+\varepsilon\ellc$, for which $\W[h_\varepsilon]=\W[h]+\varepsilon\sig>0$, and then let $\varepsilon\downarrow0$. The identity then rules out any rank $<n$.

It remains to exclude a boundary degeneration. Let $x_0\in\partial\Ct$ be such a point and diagonalize the tangential block of $\W$ there, with $e_n=\mu$. Since $\W\ge0$, all diagonal entries are nonnegative. If $\W_{nn}(x_0)=0$, the normal direction already belongs to the kernel. Otherwise $\W_{nn}(x_0)>0$, and since $\rank\W(x_0)<n$, some tangential eigenvalue must vanish, say $\W_{\alpha\alpha}(x_0)=0$. Along the inward normal geodesic the nonnegative function $\W_{\alpha\alpha}$ has a one-sided minimum at $x_0$; hence
\[
\W_{\alpha\alpha,\mu}(x_0)\le0.
\]
On the other hand, \eqref{W-normal-derivative} gives
\[
\W_{\alpha\alpha,\mu}(x_0)
=\cot\theta\,\W_{nn}(x_0)>0,
\]
a contradiction. Therefore $\W_{nn}(x_0)=0$. Moreover, for every zero tangential eigenvalue, \eqref{W-normal-derivative} then gives
\begin{equation}\label{Wnn-boundary-zero}
\W_{nn}(x_0)=0,
\qquad
\W_{\alpha\alpha,\mu}(x_0)=0.
\end{equation}
Thus the normal direction belongs to the kernel at every boundary degeneration.

Let $r_0=\rank\W(x_0)<n$ and set
\[
Q:=\sigma_{r_0+1}(\W).
\]
At $x_0$, exactly $r_0$ eigenvalues are positive. Hence $Q(x_0)=0$; moreover $Q^{ii}(x_0)=0$ in each positive eigendirection, whereas the derivatives in the zero eigendirections are equal to the product of the positive eigenvalues. By \eqref{Wnn-boundary-zero}, all zero tangential directions have vanishing normal derivative. Consequently
\begin{equation}\label{Q-normal-reduction}
Q_\mu(x_0)=Q^{nn}(x_0)\W_{nn,\mu}(x_0).
\end{equation}
To identify the last derivative, differentiate $F(\W)=R$ in the normal direction. In the positive tangential eigendirections, \eqref{W-normal-derivative} and \eqref{Wnn-boundary-zero} give
\[
\W_{ii,\mu}=-\cot\theta\,\W_{ii},
\]
while the zero tangential contributions vanish. Euler homogeneity of the one-homogeneous function $F$ therefore gives
\[
R_\mu
=F^{nn}\W_{nn,\mu}
-\cot\theta\sum_{i=1}^{r_0}F^{ii}\W_{ii}
=F^{nn}\W_{nn,\mu}-\cot\theta R.
\]
Combining this identity with \eqref{Q-normal-reduction} yields
\begin{equation*}
Q_\mu
=\frac{Q^{nn}}{F^{nn}}
\big(R_\mu+\cot\theta R\big).
\end{equation*}

The Orlicz dependence drops out of the boundary sign. Indeed, since both $h$ and $\ellc$ satisfy the same Robin condition,
\[
(\ellc/h)_\mu=0.
\]
Differentiating the definition of $R$ and using $h_\mu/h=\cot\theta$ gives
\begin{equation*}
R_\mu+\cot\theta R
=\frac{R}{kf}
\left[f_\mu+(k-1)\cot\theta f\right].
\end{equation*}
Since $H=f\ellc^{k-1}$,
\[
H_\mu=\ellc^{k-1}
\left[f_\mu+(k-1)\cot\theta f\right].
\]
Thus \eqref{OC-boundary-intro} implies $Q_\mu(x_0)\ge0$. On the other hand, $Q\ge0$ in the domain and $Q(x_0)=0$, so the one-sided outward derivative satisfies $Q_\mu(x_0)\le0$. Hence
\begin{equation}\label{Qmu-zero-final}
Q_\mu(x_0)=0.
\end{equation}
Let $\gamma(t)=\exp_{x_0}(-t\mu)$ be the inward normal geodesic. Smoothness and \eqref{Qmu-zero-final} give
\begin{equation*}
Q(\gamma(t))=O(t^2).
\end{equation*}
After shrinking the neighborhood, the $r_0$ positive eigenvalues of $\W$ stay uniformly bounded below. Since $\W\ge0$, the definition of $Q$ then implies
\begin{equation}\label{lambda-Q-final}
0\le\lambda_{\min}(\W)(\gamma(t))
\le C Q(\gamma(t))=O(t^2).
\end{equation}

We now record the boundary-uniform form of the local viscosity estimate needed below. Only the case $m=1$ of \cite[Theorem~3.4]{BISCRT} is used. Inspecting its proof, the constant in the pointwise inequality is produced from the ellipticity constants, bounds for the coefficients and their $x$-derivatives, the ambient curvature, and bounds for the fixed tensor and its first derivative; no cutoff function or positive distance from the boundary enters the calculation. Fix a closed collar $U=\{0\le d\le d_0\}$ of $x_0$. For the fixed smooth solution here, $\W$ and $\nabla\W$ are bounded on $U$, while the equation, Lemma~\ref{lem:Gamma-compact}, and smoothness of $R$ give uniform ellipticity and uniform bounds for the derivatives of $\mathscr F(A,x)=F(A)-R(x)$ required in that calculation. The curvature of the spherical cap is fixed. Thus the pointwise proof of \cite[Theorem~3.4]{BISCRT} applies at every interior point of $U$ with one constant independent of its distance to $\partial\Ct$. Consequently every smooth lower support $\xi$ of $\lambda_{\min}(\W)$ at an interior point of $U$ satisfies
\begin{equation}\label{lambda-viscosity-final}
F^{ij}\xi_{ij}
\le C\big(\xi+|\nabla\xi|\big).
\end{equation}
Because the interior degeneration has already been excluded, $\lambda_{\min}(\W)>0$ away from the boundary. Choose a geodesic ball tangent internally to $\partial\Ct$ at $x_0$ and, in the resulting annulus, use the standard barrier
\[
\psi(x)=e^{-A d(x,y)^2}-e^{-A R_0^2}.
\]
Uniform ellipticity permits $A$ to be chosen so large that
\[
F^{ij}\psi_{ij}>C(\psi+|\nabla\psi|).
\]
After multiplying by a sufficiently small positive constant, $\lambda_{\min}\ge\varepsilon\psi$ on the inner boundary and equality holds at the tangency point on the outer boundary. If $\lambda_{\min}-\varepsilon\psi$ had a negative interior minimum, a smooth translate of $\varepsilon\psi$ would be a lower support contradicting \eqref{lambda-viscosity-final}. Hence
\[
\lambda_{\min}(\W)\ge\varepsilon\psi.
\]
Along $\gamma(t)$ one has $\psi(\gamma(t))\ge c t$ for small $t$, and therefore
\begin{equation*}
\lambda_{\min}(\W)(\gamma(t))\ge c t.
\end{equation*}
This contradicts \eqref{lambda-Q-final}. Boundary degeneration is impossible, completing the full-rank argument. We obtain the following statement.

\begin{proposition}[Capillary Orlicz full rank]\label{prop:fullrank}
Let $h$ be a positive smooth $k$-admissible solution with $\W[h]\ge0$. If
\[
\mathcal A_\varphi(t)>0,
\qquad
\mathscr S_{\varphi,f}(t)\ge0
\]
for $t$ in the range of $\ellc/h$, and if $H_\mu\ge0$ on $\partial\Ct$, then $\W[h]>0$ on $\overline\Ct$.
\end{proposition}

\subsection{Orlicz-function continuation and strict convexity}

It remains to verify the structural tensor along a path beginning at a power-law equation. Put
\[
p_0:=k+1+\delta
\]
and define
\begin{equation*}
\varphi_s(t):=t^{p_0(1-s)}\varphi(t)^s,
\qquad 0\le s\le1.
\end{equation*}
Then
\[
a_s(t)=(1-s)p_0+s a_\varphi(t)\ge k+1+\delta,
\]
so the admissible solvability established above gives a unique positive admissible solution $h_s$ for every $s$. By Proposition~\ref{prop:isomorphism}, the implicit function theorem gives a local smooth solution branch in $s$; uniqueness makes these local branches agree on overlaps. Hence $s\mapsto h_s$ is continuous (indeed smooth locally) in $C^{2,\alpha}(\overline\Ct)$ along the whole interval $[0,1]$.

Assume $t a_\varphi'(t)\ge0$. Then
\[
t a_s'(t)=s\,t a_\varphi'(t)\ge0.
\]
For a general weight define
\begin{equation*}
c_\varphi(t)
:=\frac{(a_\varphi(t)-1)^2}
{a_\varphi(t)-1+t a_\varphi'(t)},
\qquad
\alpha_\varphi(t):=\frac1{k+c_\varphi(t)}.
\end{equation*}
A direct tensor calculation gives the factorization
\begin{equation}\label{S-factor-final}
\begin{aligned}
\mathscr S_{\varphi,f}(t)
={}&(k+c_\varphi(t))
\frac{
\nabla^2(\ellc H^{-\alpha_\varphi(t)})
+\ellc H^{-\alpha_\varphi(t)}\sig
}{\ellc H^{-\alpha_\varphi(t)}}\\
&+\frac{(a_\varphi-1)t a_\varphi'}
{\ellc\,[a_\varphi-1+t a_\varphi']}
\sig.
\end{aligned}
\end{equation}
Because $t a_s'\ge0$, the second term is nonnegative and
\[
0<\alpha_s(t)<\frac1k.
\]
Therefore the uniform condition \eqref{OC-intro} implies
\[
\mathscr S_{\varphi_s,f}(t)\ge0
\]
for all $s$ and all relevant $t$. Also $\mathcal A_{\varphi_s}>0$.

At $s=0$ the equation is the power problem
\[
\bk(\W[h])=\frac{f}{\ellc^{p_0}}h^{p_0-1}.
\]
For $\alpha=1/(p_0+k-1)$,
\[
\ellc H^{-\alpha}
=\left(\frac{f}{\ellc^{p_0}}\right)^{-1/(p_0+k-1)},
\]
so \eqref{OC-intro} is precisely the spherical-convexity condition used in the supercritical power full-rank theorem of \cite{MWWLp}; likewise \eqref{OC-boundary-intro} becomes its boundary structural inequality. Hence $\W[h_0]>0$.

If strict convexity were lost for the first time at $s_*$, continuity would give
\[
\W[h_{s_*}]\ge0
\]
with a zero eigenvalue. Proposition \ref{prop:fullrank}, applied using \eqref{S-factor-final}, would instead give $\W[h_{s_*}]>0$, a contradiction. Thus
\[
\W[h_s]>0\qquad\text{for all }0\le s\le1.
\]
At $s=1$, together with the admissible solvability and uniqueness proved above, this completes the proof of Theorem \ref{thm:mainA}.

\subsection{A non-power class}

The structural assumptions on the Orlicz weight are genuinely broader than the power family. For example,
\begin{equation*}
\varphi(t)=t^p(1+t^q)^\beta,
\qquad p>k+1,\quad q>0,\quad\beta>0,
\end{equation*}
satisfies
\[
a_\varphi(t)
=p+\beta q\frac{t^q}{1+t^q}>k+1,
\]
and
\[
t a_\varphi'(t)
=\beta q^2\frac{t^q}{(1+t^q)^2}>0.
\]
Thus Theorem \ref{thm:mainA} applies to genuinely non-power Orlicz weights whenever the prescribed datum satisfies the capillary Orlicz-convexity condition.

\begin{remark}[Endpoint and open directions]
The argument is intentionally restricted to $0<\theta<\pi/2$. At $\theta=\pi/2$ one needs a separate Neumann/reflection regularity treatment, while for $\theta>\pi/2$ the sign used in the boundary $C^2$ barrier is lost. A full intermediate-order capillary Orlicz--Brunn--Minkowski inequality and the obtuse-angle problem remain natural geometric extensions of the present elliptic theory.
\end{remark}

\medskip
\noindent\textbf{Acknowledgments:}
Hua Zhu was supported by the National Natural Science Foundation of China (Grant No. 12501273, 12661041) and the Research Foundation of Southwest University of Science and Technology (Grant No. 25zx7153).

\medskip
\noindent\textbf{AI Statement:}
During the development and preparation of this work, the author used a generative AI tool for preliminary, non-authoritative exploratory assistance, including organizational discussion and consideration of possible approaches. AI-generated outputs were not treated as mathematical sources, and no argument was included without independent verification and substantial revision by the author. All theorem statements, proofs, computations, references, and final arguments were independently checked and finalized by the author, who takes full responsibility for the correctness, originality, and integrity of the work.


\begin{thebibliography}{99}

\bibitem{BISCRT}
P. Bryan, M. N. Ivaki, and J. Scheuer,
\emph{Constant rank theorems for curvature problems via a viscosity approach},
Calc. Var. Partial Differential Equations \textbf{62} (2023), Article 98.

\bibitem{CNS}
L. Caffarelli, L. Nirenberg, and J. Spruck,
\emph{The Dirichlet problem for nonlinear second-order elliptic equations III: Functions of the eigenvalues of the Hessian},
Acta Math. \textbf{155} (1985), 261--301.

\bibitem{GHW}
R. J. Gardner, D. Hug, and W. Weil,
\emph{The Orlicz--Brunn--Minkowski theory: a general framework, additions, and inequalities},
J. Differential Geom. \textbf{97} (2014), 427--476.

\bibitem{GuanMa}
P. Guan and X.-N. Ma,
\emph{The Christoffel--Minkowski problem I: Convexity of solutions of a Hessian equation},
Invent. Math. \textbf{151} (2003), 553--577.

\bibitem{GuanMaZhou}
P. Guan, X.-N. Ma, and F. Zhou,
\emph{The Christoffel--Minkowski problem III: Existence and convexity of admissible solutions},
Comm. Pure Appl. Math. \textbf{59} (2006), 1352--1376.

\bibitem{HLYZ}
C. Haberl, E. Lutwak, D. Yang, and G. Zhang,
\emph{The even Orlicz Minkowski problem},
Adv. Math. \textbf{224} (2010), 2485--2510.

\bibitem{HuIvakiLpCM}
Y. Hu and M. N. Ivaki,
\emph{Capillary $L_p$-Christoffel--Minkowski problem},
arXiv:2512.15464v2, 2025.

\bibitem{HuIvaki}
Y. Hu and M. N. Ivaki,
\emph{Capillary $L_p$-curvature problem},
arXiv:2602.21832, 2026.

\bibitem{HuIvakiScheuer}
Y. Hu, M. N. Ivaki, and J. Scheuer,
\emph{Capillary Christoffel--Minkowski problem},
arXiv:2504.09320, 2025.

\bibitem{LiebermanTrudinger}
G. M. Lieberman and N. S. Trudinger,
\emph{Nonlinear oblique boundary value problems for nonlinear elliptic equations},
Trans. Amer. Math. Soc. \textbf{295} (1986), 509--546.

\bibitem{Lutwak93}
E. Lutwak,
\emph{The Brunn--Minkowski--Firey theory I: Mixed volumes and the Minkowski problem},
J. Differential Geom. \textbf{38} (1993), 131--150.

\bibitem{CMCap}
X. Mei, G. Wang, and L. Weng,
\emph{The capillary Christoffel--Minkowski problem},
Calc. Var. Partial Differential Equations \textbf{65} (2026), Article 186.

\bibitem{MinkowskiCap}
X. Mei, G. Wang, and L. Weng,
\emph{The capillary Minkowski problem},
Adv. Math. \textbf{469} (2025), Article 110230.

\bibitem{MWWLp}
X. Mei, G. Wang, and L. Weng,
\emph{Prescribed $L_p$ curvature problem for convex capillary hypersurface},
arXiv:2512.16686v3, revised 27 September 2026.

\bibitem{MWWLpMinkowski}
X. Mei, G. Wang, and L. Weng,
\emph{The capillary $L_p$-Minkowski problem},
arXiv:2505.07746v3, revised 22 June 2026.

\bibitem{AFII}
X. Mei, G. Wang, L. Weng, and C. Xia,
\emph{Alexandrov--Fenchel inequalities for convex hypersurfaces in the half-space with capillary boundary II},
Math. Z. \textbf{310} (2025), Article 71.

\bibitem{WangWengXia}
G. Wang, L. Weng, and C. Xia,
\emph{A Minkowski-type inequality for capillary hypersurfaces in a half-space},
J. Funct. Anal. \textbf{287} (2024), Article 110496.

\bibitem{WangZhu}
X. Wang and B. Zhu,
\emph{The capillary Orlicz--Minkowski problem},
J. Funct. Anal. \textbf{291} (2026), Article 111574.

\bibitem{XiJinLeng}
D. Xi, H. Jin, and G. Leng,
\emph{The Orlicz Brunn--Minkowski inequality},
Adv. Math. \textbf{260} (2014), 350--374.
\end{thebibliography}
\end{document}